\documentstyle{article}
\newtheorem{rem}{Remark}
\newtheorem{theor}{Theorem}
\newtheorem{dfn}{Definition}
\begin{document}
\begin{center}
\Large{\bf  ON THE THEORY OF SPACES OF CONSTANT \\[2mm] CURVATURE
}\vspace{4mm}\normalsize
\end{center}
 \begin{center}
\Large{\bf Valery Dryuma}\vspace{4mm}\normalsize
\end{center}
\begin{center}
{\bf Institute of Mathematics and Informatics AS Moldova, Kishinev}\vspace{4mm}\normalsize
\end{center}
\begin{center}
{\bf  E-mail: valery@dryuma.com;\quad cainar@mail.md}\vspace{4mm}\normalsize
\end{center}
\begin{center}
{\bf  Abstract}\vspace{4mm}\normalsize
\end{center}

   Some examples of three-dimensional metrics of constant curvature
   defined by the solutions of nonlinear integrable differential equations and
   their generalizations are constructed. The properties of  Riemann
   extensions of the metrics of constant curvature are studied.
   The connection with the theory of normal Riemann spaces are discussed.

\section{Introduction}

   The metrics of constant curvature are defined
by a following condition on a curvature tensor
\begin{equation} \label{dryuma:eq1}
R_{i j k l}-\lambda\left(g_{i k}g_{j l}-g_{i l}g_{j k}\right)=0.
\end{equation}

   In the case of diagonal metric
\[
ds^2= A(x,y,z)^2dx^2+B(x,y,z)^2dy^2+C(x,y,z)^2dz^2
\]
the conditions  (\ref{dryuma:eq1}) are equivalent the  system of
differential equations on the functions $A,B,C$.

   They looks as
\begin{equation} \label{dryuma:eq2}
{\frac {\partial ^{2}}{\partial y\partial z}}A(x,y,z)={\frac
{\left ({ \frac {\partial }{\partial z}}B(x,y,z)\right ){\frac
{\partial }{
\partial y}}A(x,y,z)}{B(x,y,z)}}+{\frac {\left ({\frac {\partial }{
\partial y}}C(x,y,z)\right ){\frac {\partial }{\partial z}}A(x,y,z)}{C
(x,y,z)}} ,
\end{equation}
\begin{equation} \label{dryuma:eq3}
{\frac {\partial ^{2}}{\partial x\partial z}}B(x,y,z)={\frac
{\left ({ \frac {\partial }{\partial z}}A(x,y,z)\right ){\frac
{\partial }{
\partial x}}B(x,y,z)}{A(x,y,z)}}+{\frac {\left ({\frac {\partial }{
\partial x}}C(x,y,z)\right ){\frac {\partial }{\partial z}}B(x,y,z)}{C
(x,y,z)}} ,
\end{equation}
\begin{equation} \label{dryuma:eq4}
{\frac {\partial ^{2}}{\partial x\partial y}}C(x,y,z)={\frac
{\left ({ \frac {\partial }{\partial y}}A(x,y,z)\right ){\frac
{\partial }{
\partial x}}C(x,y,z)}{A(x,y,z)}}+{\frac {\left ({\frac {\partial }{
\partial x}}B(x,y,z)\right ){\frac {\partial }{\partial y}}C(x,y,z)}{B
(x,y,z)}}
\end{equation}
\begin{equation} \label{dryuma:eq5}
\lambda\,C(x,y,z)B(x,y,z)+{\frac {\left ({\frac {\partial
}{\partial x}}C(x,y,z)\right ){\frac {\partial }{\partial
x}}B(x,y,z)}{\left (A(x ,y,z)\right )^{2}}}+{\frac {{\frac
{\partial ^{2}}{\partial {z}^{2}}}B
(x,y,z)}{C(x,y,z)}}-\]\[-{\frac {\left ({\frac {\partial
}{\partial z}}B(x, y,z)\right ){\frac {\partial }{\partial
z}}C(x,y,z)}{\left (C(x,y,z) \right )^{2}}}+{\frac {{\frac
{\partial ^{2}}{\partial {y}^{2}}}C(x,y, z)}{B(x,y,z)}}-{\frac
{\left ({\frac {\partial }{\partial y}}B(x,y,z) \right ){\frac
{\partial }{\partial y}}C(x,y,z)}{\left (B(x,y,z) \right )^{2}}}
=0,
\end{equation}
\\[1mm]
\begin{equation} \label{dryuma:eq6}
\lambda\,A(x,y,z)C(x,y,z)-{\frac {\left ({\frac {\partial
}{\partial z}}A(x,y,z)\right ){\frac {\partial }{\partial
z}}C(x,y,z)}{\left (C(x ,y,z)\right )^{2}}}+{\frac {{\frac
{\partial ^{2}}{\partial {z}^{2}}}A
(x,y,z)}{C(x,y,z)}}+\]\[+{\frac {\left ({\frac {\partial
}{\partial y}}A(x, y,z)\right ){\frac {\partial }{\partial
y}}C(x,y,z)}{\left (B(x,y,z) \right )^{2}}}+{\frac {{\frac
{\partial ^{2}}{\partial {x}^{2}}}C(x,y, z)}{A(x,y,z)}}-{\frac
{\left ({\frac {\partial }{\partial x}}A(x,y,z) \right ){\frac
{\partial }{\partial x}}C(x,y,z)}{\left (A(x,y,z) \right
)^{2}}}=0,
\end{equation}
\\[1mm]
\begin{equation} \label{dryuma:eq7}
\lambda\,B(x,y,z)A(x,y,z)-{\frac {\left ({\frac {\partial
}{\partial x}}A(x,y,z)\right ){\frac {\partial }{\partial
x}}B(x,y,z)}{\left (A(x ,y,z)\right )^{2}}}+{\frac {{\frac
{\partial ^{2}}{\partial {x}^{2}}}B
(x,y,z)}{A(x,y,z)}}+\]\[+{\frac {\left ({\frac {\partial
}{\partial z}}A(x, y,z)\right ){\frac {\partial }{\partial
z}}B(x,y,z)}{\left (C(x,y,z) \right )^{2}}}+{\frac {{\frac
{\partial ^{2}}{\partial {y}^{2}}}A(x,y, z)}{B(x,y,z)}}-{\frac
{\left ({\frac {\partial }{\partial y}}A(x,y,z) \right ){\frac
{\partial }{\partial y}}B(x,y,z)}{\left (B(x,y,z) \right
)^{2}}}=0.
\end{equation}

     The solutions of full system of equations (\ref{dryuma:eq2}-\ref{dryuma:eq7})
     depends from the values of $\lambda$ and are important in geometry of 3-dimensional
     spaces.

     Remark that first three equations do not depend from the parameter $\lambda$
     and their solutions are common for all possible cases of the equations
     described the metrics of the spaces of constant curvature.

    The system of equations (\ref{dryuma:eq2}-\ref{dryuma:eq4})
    as example of exactly integrable system of multidimensional equations was
    discovered in the work of author (\cite{dryuma:darboux}) and have been studied in
     (\cite{dryuma:darboux1}-\cite{dryuma:darboux2}) together with their applications at
     the problem of three-orthogonal system of surfaces in a $RP^3$ and $R^3$ spaces.

\section{The space of zero curvature $\lambda=0$}

   The properties of the spaces of curvature $\lambda=0$ in diagonal metrics have
   been studied in (\cite{dryuma:zakh}).

   Here we consider the metric in non diagonal form (\cite{dryuma:wolf})
\begin{equation} \label{dryuma:eq8}
{{\it ds}}^{2}={y}^{2}{{\it dx}}^{2}+2\,\left
(l(x,z){y}^{2}+m(x,z) \right ){\it dx}\,{\it dz}+2\,{\it dy}\,{\it
dz}+\]\[+\left (\left (l(x,z) \right )^{2}{y}^{2}-2\,\left ({\frac
{\partial }{\partial x}}l(x,z) \right
)y+2\,l(x,z)m(x,z)+2\,l(x,z)\right ){{\it dz}}^{2}
\end{equation}
 with some functions $l(x,z)$ and $m(x,z)$.

   The condition on the curvature tensor
\[
R_{i j k l}=0
\]
for the metric (\ref{dryuma:eq8}) lead to the relations
\[
R_{1313}=\left ({\frac {\partial ^{3}}{\partial
{x}^{3}}}l(x,z)-3\,l(x,z){ \frac {\partial }{\partial
x}}l(x,z)+{\frac {\partial }{\partial z}}l( x,z)\right
){y}+\]\[+\left (\!-\!l(x,z){\frac {\partial ^{2}}{\partial {x}^
{2}}}m(x,z)\!+\!{\frac {\partial ^{2}}{\partial x\partial
z}}m(x,z)\!-\!3\, \left ({\frac {\partial }{\partial
x}}l(x,z)\right ){\frac {\partial } {\partial
x}}m(x,z)\!-\!2\,m(x,z){\frac {\partial ^{2}}{\partial {x}^{2}}}
l(x,z)\!-\!{\frac {\partial ^{2}}{\partial {x}^{2}}}l(x,z)\right
)\!-\!\]\[\left(-m(x,z) {\frac {\partial }{\partial
z}}m(x,z)\!+\!m(x,z)l(x,z){\frac {\partial }{
\partial x}}m(x,z)\!+\!\left ({\frac {\partial }{\partial x}}l(x,z)\right
)m(x,z)\!+\!2\,\left ({\frac {\partial }{\partial x}}l(x,z)\right
)\left ( m(x,z)\right )^{2}\right)/y =0,
\]
and
\[
R_{1323}=\left({\frac {\partial }{\partial x}}l(x,z)-{\frac
{\partial }{\partial z}}m (x,z)+2\,\left ({\frac {\partial
}{\partial x}}l(x,z)\right )m(x,z)+l( x,z){\frac {\partial
}{\partial x}}m(x,z)\right)/y=0
\]
which are equivalent the system of equations for the functions
$l(x,z)$ and $m(x,z)$
\begin{equation} \label{dryuma:eq9}
{\frac {\partial ^{3}}{\partial {x}^{3}}}l(x,z)-3\,l(x,z){\frac {
\partial }{\partial x}}l(x,z)+{\frac {\partial }{\partial z}}l(x,z)
=0,
\end{equation}
\begin{equation} \label{dryuma:eq10}
{\frac {\partial }{\partial x}}l(x,z)-{\frac {\partial }{\partial
z}}m (x,z)+2\,\left ({\frac {\partial }{\partial x}}l(x,z)\right
)m(x,z)+l( x,z){\frac {\partial }{\partial x}}m(x,z) =0.
\end{equation}

     So we can formulate the following
\begin{theor}
  There is exists a class of 3-dimensional flat  metrics  defined by the
  solutions of the system of equations (\ref{dryuma:eq9}-\ref{dryuma:eq10}).
\end{theor}

 Remark that the first equation of the system  is the famous
     KdV-equation and this fact may be used for studying of the
     properties of orthogonal metrics.

    Let us consider some examples.

    1.
\[
l(x,z)=-\frac{x}{3z},\quad m(x,z)=-1/2+{\it F_1}({\frac {z}{{x}^{3}}}){x}^{-2}.
\]

    2.
\[
l(x,z)=-4\,\left (\cosh(x-4\,z)\right )^{-2} , \quad
m(x,z)=-\frac{1}{2}
\]
\begin{rem}

     In the simplest case three-dimensional metrics of zero curvature look as
\[
ds^2={y}^{2}{{\it dx}}^{2}+2\,\left (l(x,z){y}^{2}-1/2\right ){\it
dx}\,{ \it dz}+2\,{\it dy}\,{\it dz}+\]\[+\left (\left
(l(x,z)\right )^{2}{y}^{2}- 2\,\left ({\frac {\partial }{\partial
x}}l(x,z)\right )y+l(x,z)\right ){{\it dz}}^{2}
\]
where the function $l(x,z)$ is solution of classical KdV-equation
\[
{\frac {\partial }{\partial z}}l(x,z)-3\,l(x,z){\frac {
\partial }{\partial x}}l(x,z)+{\frac {\partial ^{3}}{\partial {x}^{3}}
}l(x,z)=0.
\]

   In particular the functions
\[
l(x,z)=-4\,\left (\cosh(x-4\,z)\right )^{-2},
\]
and
\[
l(x,z)=-24\,{\frac {4\,\cosh(2\,x-8\,z)+\cosh(4\,x-64\,z)+3}{\left
(3\,\cosh( x-28\,z)+\cosh(3\,x-36\,z)\right )^{2}}}
\]
give us the examples of such type of metrics.

     In spite of the fact that the determinant of the metric do not depends from
     the function $l(x,z)$ it is possible to distinguish the properties of metrics
     with the help of eigenvalue equation for the Laplace operator defined on the
     1-forms
\[
A(x,y,z)=A_i(x,y,z)dx^i.
\]

   It has the form
\[
g^{ij}\nabla_i\nabla_j A_k-R^l_k A_l=-\lambda A_k.
\]

    In particular case $l(x,z)=0$ and $A_i(x,y,z)=[h(y),q(y),f(y)]$ these equations
    take the form
\[
-1/4\,{\frac {{\frac {d}{dy}}h(y)-4\,\left ({\frac
{d}{dy}}f(y)\right ){y}^{2}-\left ({\frac
{d^{2}}{d{y}^{2}}}h(y)\right )y+4\,\lambda\,h(y
){y}^{3}}{{y}^{3}}}=0,
\]
\[
-1/4\,{\frac {-6\,h(y)-3\,q(y)+3\,\left ({\frac {d}{dy}}q(y)\right
)y+ 4\,\left ({\frac {d}{dy}}h(y)\right )y-\left ({\frac
{d^{2}}{d{y}^{2}} }q(y)\right
){y}^{2}+4\,f(y){y}^{2}+4\,\lambda\,q(y){y}^{4}}{{y}^{4}}}=0,
\]
\[
-1/4\,{\frac {{\frac {d}{dy}}f(y)-\left ({\frac
{d^{2}}{d{y}^{2}}}f(y) \right
)y+4\,\lambda\,f(y){y}^{3}}{{y}^{3}}}=0.
\]

    The solutions of this system are dependent from the eigenvalues $\lambda$ and
    characterize the properties of corresponding flat metrics.
\end{rem}

\begin{rem}

 In theory of varieties the Chern-Simons characteristic class is constructed from a
 matrix gauge connection $A^i_{jk}$ as
\[
W(A)=\frac{1}{4\pi^2}\int d^3x
\epsilon^{ijk}tr\left(\frac{1}{2}A_i\partial_j
A_k+\frac{1}{3}A_iA_jA_k \right).
\]

   This term can be translated into three-dimensional geometric quantity by replacing the matrix connection
   $A^i_{jk}$ with the Christoffel connection $\Gamma^i_{jk}$.

   For the density of  Chern-Simons invariant can be obtained the
   expression (\cite{dryuma:jac})
\begin{equation} \label{dryuma:eq11}
CS(\Gamma)=\epsilon^{i j k}(\Gamma^p_{i q}\Gamma^q_{k
p;j}+\frac{2}{3}\Gamma^p_{i q}\Gamma^q_{j r}\Gamma^r_{k p}).
\end{equation}

   For the metric (\ref{dryuma:eq8}) we find the quantity
\[
CS(\Gamma)=-{\frac {5\,l(x,z){\frac {\partial }{\partial
x}}l(x,z)-5\,{\frac {
\partial }{\partial z}}l(x,z)}{\sqrt {{y}^{2}}}}-\]\[-{\frac {3\,l(x,z){
\frac {\partial }{\partial x}}m(x,z)-4\,\left ({\frac {\partial }{
\partial x}}l(x,z)\right )m(x,z)-3\,{\frac {\partial }{\partial z}}m(x
,z)-2\,{\frac {\partial }{\partial x}}l(x,z)}{{y}^{2}\sqrt
{{y}^{2}}}}.
\]

     Using this formulae for  the first example one get
\[
CS(\Gamma)={\frac {10}{9}}\,{\frac {x{\it
csgn}(y)}{{z}^{2}y}}-10/3\,{\it \_F1}({ \frac {z}{{x}^{3}}}){\it
csgn}(y){z}^{-1}{x}^{-2}{y}^{-3}.
\]

    For the second  example this quantity is
\[
CS(\Gamma)=160\,{\frac {\cosh(x-4\,z)\left (\sinh(x-4\,z)\right
)^{3}{\it csgn}(y )}{y}}.
\]
\end{rem}

\section{The metrics of nonzero  curvature $\lambda\neq0$}

     The metric of the space of positive curvature $\lambda=1$
 is defined by
\begin{equation} \label{dryuma:eq12}
ds^2=\frac{dx^2+dy^2+dz^2}{(1+(x^2+y^2+z^2)/4))^2}.
\end{equation}

     The metric of the space of negative curvature $\lambda=-1$
 is defined by
\begin{equation} \label{dryuma:eq13}
ds^2=\frac{dx^2+dy^2+dz^2}{z^2}.
\end{equation}

     Starting from these expressions it is possible to get
     more general examples of the metrics of constant curvature.

     For example, the substitution
\[
A(x,y,z)=\frac{1}{z},\quad  B(x,y,z)=\frac{1}{z}+v(z),\quad
C(x,y,z)=\frac{1}{z}+u(z)
\]
into the system (\ref{dryuma:eq2})-(\ref{dryuma:eq7}) and
integration of corresponding  equations for the functions $u(z)$
and $v(z)$ give us the metric
\[
ds^2={\frac {{{\it dx}}^{2}}{{z}^{2}}}+\left ({z}^{-2}-1\right
){{\it dy}}^ {2}+{\frac {{{\it dz}}^{2}}{{z}^{2}\left
(1-{z}^{2}\right )}}
\]
of negative curvature $\lambda=-1$.

    The substitution
\[
A(x,y,z)=\left (1+1/4\,{x}^{2}+1/4\,{y}^{2}+1/4\,{z}^{2}\right
)^{-1}+ U(x,y,z)
\]
into the system ~(\ref{dryuma:eq2})-(\ref{dryuma:eq7}) and integration of
corresponding equations lead to the metrics
\[
ds^2=A(x,y,z)^2dx^2+\frac{dy^2+dz^2}{(1+(x^2+y^2+z^2)/4))^2},
\]
with the function $A(x,y,z)$ in form
\[
A(x,y,z)=4\,{\frac {\sqrt {{x}^{2}\left
({x}^{2}-4\,\lambda+4\right )}}{\left ( {x}^{2}-4\,\lambda+4\right
)\left (4+{x}^{2}+{y}^{2}+{z}^{2}\right )}}
\]

  The Chern-Simons invariant of orthogonal 3-dimensional  metrics with conditions (\ref{dryuma:eq2})-(\ref{dryuma:eq4})
  on components is given  by the expression
\begin{equation} \label{dryuma:eq14}
CS(\Gamma)=10\,{\frac {\left ({\frac {\partial }{\partial
y}}A(x,y,z)\right ) \left ({\frac {\partial }{\partial
z}}B(x,y,z)\right ){\frac {
\partial }{\partial x}}C(x,y,z)-\left ({\frac {\partial }{\partial z}}
A(x,y,z)\right )\left ({\frac {\partial }{\partial
x}}B(x,y,z)\right ) {\frac {\partial }{\partial
y}}C(x,y,z)}{A(x,y,z)B(x,y,z)C(x,y,z)}}
\end{equation}

   After substitution of the components of metric
   (\ref{dryuma:eq13}) into this expression  one obtain
\[
CS(\Gamma)=0.
\]

     As it is follows from (\ref{dryuma:eq14}) the class of 3-dimensional orthogonal metrics
     with vanishing Chern-Simons invariant is defined by the condition
 \begin{equation} \label{dryuma:eq15}
  \left( {\frac {\partial }{\partial z}}A \left( x,y,z \right)
 \right)  \left( {\frac {\partial }{\partial x}}B \left( x,y,z
 \right)  \right) {\frac {\partial }{\partial y}}C \left( x,y,z
 \right) - \left( {\frac {\partial }{\partial y}}A \left( x,y,z
 \right)  \right)  \left( {\frac {\partial }{\partial z}}B \left( x,y,
z \right)  \right) {\frac {\partial }{\partial x}}C \left( x,y,z
 \right)=0.
 \end{equation}

    The metrics of the form
\[
ds^2 =2 E(x,y,z)dxdy+dz^2
\]
are the metrics of constant curvature if the function $E(x,y,z)$ is defined by
\[
E(x,y,z)=1/4\,{\frac {\left (F(x,y)\sin(\sqrt
{\lambda}z)-F(x,y)\cos( \sqrt {\lambda}z)\right )^{2}}{\lambda}},
\]
where
\[
-4\,\left ({\frac {\partial }{\partial y}}F(x,y)\right ){\frac {
\partial }{\partial x}}F(x,y)+4\,F(x,y){\frac {\partial ^{2}}{
\partial x\partial y}}F(x,y)+\left (F(x,y)\right )^{4}=0,
\]
or
\[
4\,{\frac {\partial ^{2}}{\partial x\partial
y}}U(x,y)+{e^{2\,U(x,y)}} =0
\]
after the substitution
\[
F(x,y)=\exp(U(x,y)).
\]

 Last equation is the famous Liouville equation with general
 solution
 \[
U(x,y)=1/2\,\ln (-4\,{\frac {\left ({\frac {d}{dx}}a(x)\right
){\frac {d}{dy}}b(y)}{\left (\left (a(x)\right )^{2}+\left
(b(y)\right )^{2} \right )^{2}}}).
 \]

   The Chern-Simons invariant for this example of metrics is
\[
CS(\Gamma)=0.
\]

    Let us consider the metrics of the form
\[
ds^2=dx^2+2\cos(u(x,y)dxdy+A(x,y)^2dz^2.
\]

   The conditions (\ref{dryuma:eq1}) lead to the linear system of
   equations
\[
 {\frac {\partial ^{2}}{\partial x\partial y}}A(x,y)+\lambda\,A(x,y)
\cos(u(x,y))=0,
\]
\[
{\frac {\partial ^{2}}{\partial {x}^{2}}}A(x,y)-{\frac
{\cos(u(x,y)) \left ({\frac {\partial }{\partial x}}u(x,y)\right
){\frac {\partial } {\partial
x}}A(x,y)}{\sin(u(x,y))}}+\lambda\,A(x,y)+{\frac {\left ({ \frac
{\partial }{\partial x}}u(x,y)\right ){\frac {\partial }{
\partial y}}A(x,y)}{\sin(u(x,y))}}=0,
\]
\[
{\frac {\partial ^{2}}{\partial {y}^{2}}}A(x,y)-{\frac
{\cos(u(x,y)) \left ({\frac {\partial }{\partial y}}u(x,y)\right
){\frac {\partial } {\partial
y}}A(x,y)}{\sin(u(x,y))}}+\lambda\,A(x,y)+{\frac {\left ({ \frac
{\partial }{\partial x}}A(x,y)\right ){\frac {\partial }{
\partial y}}u(x,y)}{\sin(u(x,y))}}=0,
\]
which is compatible on the solutions of the $"\sin-Gordon"$ equation
 $(\lambda=-\mu)$
\[
{\frac {\partial ^{2}}{\partial x\partial
y}}u(x,y)+\lambda\,\sin(u(x, y))=0
\]

    So for any solution of the $"\sin-Gordon"$ equation  one possible to find the
    function $A(x,y)$ with the help of solution of the corresponding linear system
    of equations.

    To take one example.

    The simplest solution of the equation

\[
{\frac {\partial ^{2}}{\partial x\partial y}}u(x,y)-\sin(u(x,
y))=0
\]
is given by
\[
u(x,y)=4\arctan(\exp(x+y)).
\]

    By this condition the linear system looks as
\[
{\frac {\partial ^{2}}{\partial x\partial y}}A(x,y)+{\frac {\left
(6\, {e^{2\,x+2\,y}}-1-{e^{4\,x+4\,y}}\right )A(x,y)}{\left
(1+{e^{2\,x+2\, y}}\right )^{2}}}=0,
\]
\[
\left (1\!-\!{e^{4\,x+4\,y}}\right )A(x,y)\!+\!\left
(-2\,{e^{2\,x+2\,y}}-{e^{ 4\,x+4\,y}}-1\right ){\frac {\partial
}{\partial y}}A(x,y)\!+\!\left (-6\,
{e^{2\,x+2\,y}}\!+\!1\!+\!{e^{4\,x+4\,y}}\right ){\frac {\partial
}{\partial x }}A(x,y)+\]\[+\left (-1+{e^{4\,x+4\,y}}\right ){\frac
{\partial ^{2}}{
\partial {x}^{2}}}A(x,y)=0,
\]
\[
\left (1\!-\!{e^{4\,x+4\,y}}\right )A(x,y)\!+\!\left
(-6\,{e^{2\,x+2\,y}}+1+{e ^{4\,x+4\,y}}\right ){\frac {\partial
}{\partial y}}A(x,y)+\left (-2\,
{e^{2\,x+2\,y}}-{e^{4\,x+4\,y}}-1\right ){\frac {\partial
}{\partial x }}A(x,y)+\]\[+\left (-1+{e^{4\,x+4\,y}}\right ){\frac
{\partial ^{2}}{
\partial {y}^{2}}}A(x,y)=0
\]

    Its solution is
\[
A(x,y)={\frac {{e^{x+y}}}{1+{e^{2\,x+2\,y}}}}.
\]

\section{An examples of essentially three-dimensional metrics}

  We consider a family of three-dimensional  metrics in form
\[
{{\it ds}}^{2}={{\it dx}}^{2}+2\,\cos(B(x,y,z)){\it dxdy}+{{\it
dy}}^{ 2}+\left ({\frac {\partial }{\partial z}}B(x,y,z)\right
)^{2}{{\it dz} }^{2}
\]
 The conditions (\ref{dryuma:eq1}) for a such metrics lead to a following  system of equations
\[
{\frac {\partial ^{2}}{\partial x\partial
y}}B(x,y,z)+1/4\,\sin(B(x,y, z))\left (-1+4\,\lambda\right )=0,
\]
\[
{\frac {\partial ^{3}}{\partial {y}^{2}\partial z}}B(x,y,z)-{\frac
{ \cos(B(x,y,z))\left ({\frac {\partial }{\partial
y}}B(x,y,z)\right ){ \frac {\partial ^{2}}{\partial y\partial
z}}B(x,y,z)}{\sin(B(x,y,z))}} +{\frac {\left ({\frac {\partial
^{2}}{\partial x\partial z}}B(x,y,z) \right ){\frac {\partial
}{\partial y}}B(x,y,z)}{\sin(B(x,y,z))}}-\]\[- \left
(1/4-\lambda\right ){\frac {\partial }{\partial z}}B(x,y,z)=0,
\]
\[
{\frac {\partial ^{3}}{\partial {x}^{2}\partial z}}B(x,y,z)-{\frac
{ \cos(B(x,y,z))\left ({\frac {\partial }{\partial
x}}B(x,y,z)\right ){ \frac {\partial ^{2}}{\partial x\partial
z}}B(x,y,z)}{\sin(B(x,y,z))}} +{\frac {\left ({\frac {\partial
}{\partial x}}B(x,y,z)\right ){\frac {\partial ^{2}}{\partial
y\partial z}}B(x,y,z)}{\sin(B(x,y,z))}}-\]\[- \left
(1/4-\lambda\right ){\frac {\partial }{\partial z}}B(x,y,z)=0.
\]

     This system is compatible and its solutions gives us the
     examples of three-dimensional metrics of constant
     curvature.

 The Chern-Simons term in this case is defined by
\[
CS(\Gamma)=-5/2\,\left ({\frac {\partial ^{2}}{\partial x\partial
z}}B(x,y,z) \right ){\frac {\partial }{\partial
y}}B(x,y,z)+5/2\,\left ({\frac {
\partial }{\partial x}}B(x,y,z)\right ){\frac {\partial ^{2}}{
\partial y\partial z}}B(x,y,z).
\]

     Let us consider an example.

    For the sake of simplicity we present the metric in equivalent form
\[
ds^2={d{{x}}}^{2}+2\,u(x,y,z)d{{x}}d{{y}}+{d{{y}}}^{2}+{\frac
{\left ({ \frac {\partial }{\partial z}}u(x,y,z)\right
)^{2}{d{{z}}}^{2}}{1- \left (u(x,y,z)\right )^{2}}},
\]
where the function $u(x,y,z)$ is determined from the condition
\[
u(x,y,z)=\arccos(B(x,y,z)).
\]

   A system of equations for the function $B(x,y,z)$ grade at that into more
   simple system on the function $u(x,y,z)$ .

   Its integration at the condition $\lambda=1/4$ lead to the metric
\[
{{\it ds}}^{2}={{\it dx}}^{2}+1/2\,{\frac {\left ({z}^{4}\left
(f(x) \right )^{2}\left (h(y)\right )^{2}+1\right ){\it dx}\,{\it
dy}}{{z}^{ 2}f(x)h(y)}}+{{\it dy}}^{2}-4\,{\frac {{{\it
dz}}^{2}}{{z}^{2}}}
\]
where $f(x),h(y)$ are arbitrary functions.

The Ricci tensor of this case has a form
\[
\mbox {{R}}_{{a}}\mbox {{}}_{{b}}=\left [\begin {array} {ccc}
1/2&1/4\,{\frac {{z}^{4}\left (f(x)\right )^{2}\left (h(y) \right
)^{2}+1}{{z}^{2}f(x)h(y)}}&0\\\noalign{\medskip}1/4\,{\frac {{z
}^{4}\left (f(x)\right )^{2}\left (h(y)\right
)^{2}+1}{{z}^{2}f(x)h(y)
}}&1/2&0\\\noalign{\medskip}0&0&-2\,{z}^{-2}\end {array}\right ].
\]

   The Chern-Simons invariant of this metric equal to zero.

\begin{rem}
   The absolute value of curvature $\lambda=1/4$ is special in
   theory of Riemann manifolds and is connected with the theory of manifolds with
    pinched curvature.
\end{rem}
\section{The Rieman extensions of the spaces \\[1mm] of constant
curvature $\lambda=\pm1,0$}

  The space with diagonal metrics has a following coefficients of connections
\begin{equation} \label{dryuma:eq16}
\Gamma^1_{11}={\frac {{\frac {\partial }{\partial
x}}A(x,y,z)}{A(x,y,z)}},\quad \Gamma^2_{11}=-{\frac
{A(x,y,z){\frac {\partial }{\partial y}}A(x,y,z)}{\left (B(x,y
,z)\right )^{2}}},\Gamma^3_{11}=-{\frac {A(x,y,z){\frac {\partial
}{\partial z}}A(x,y,z)}{\left (C(x,y ,z)\right )^{2}}} \]\[\quad
\Gamma^1_{12}={\frac {{\frac {\partial }{\partial
y}}A(x,y,z)}{A(x,y,z)}},\quad \Gamma^2_{12}= {\frac {{\frac
{\partial }{\partial x}}B(x,y,z)}{B(x,y,z)}} ,\quad
\Gamma^3_{12}=0,\quad \Gamma^1_{13}= {\frac {{\frac {\partial
}{\partial z}}A(x,y,z)}{A(x,y,z)}} ,\]\[\quad
\Gamma^2_{13}=0,\quad \Gamma^3_{13}= {\frac {{\frac {\partial
}{\partial x}}C(x,y,z)}{C(x,y,z)}} ,\quad \Gamma^1_{22}=-{\frac
{B(x,y,z){\frac {\partial }{\partial x}}B(x,y,z)}{\left (A(x,y
,z)\right )^{2}}} ,\quad \Gamma^2_{22}={\frac {{\frac {\partial
}{\partial y}}B(x,y,z)}{B(x,y,z)}} ,\]\[\quad
\Gamma^3_{22}=-{\frac {B(x,y,z){\frac {\partial }{\partial
z}}B(x,y,z)}{\left (C(x,y ,z)\right )^{2}}} ,\quad
\Gamma^1_{23}=0,\quad \Gamma^2_{23}={\frac {{\frac {\partial
}{\partial z}}B(x,y,z)}{B(x,y,z)}} ,\quad \Gamma^3_{23}={\frac
{{\frac {\partial }{\partial y}}C(x,y,z)}{C(x,y,z)}},\]\[\quad
\Gamma^1_{33}=-{\frac {C(x,y,z){\frac {\partial }{\partial
x}}C(x,y,z)}{\left (A(x,y ,z)\right )^{2}}}, \quad
\Gamma^2_{33}=-{\frac {C(x,y,z){\frac {\partial }{\partial
y}}C(x,y,z)}{\left (B(x,y ,z)\right )^{2}}} \quad
\Gamma^3_{33}={\frac {{\frac {\partial }{\partial
z}}C(x,y,z)}{C(x,y,z)}}.
\end{equation}

   Using these expressions we introduce a six-dimensional metric
\begin{equation} \label{dryuma:eq17}
^{6}ds^2=-2\Gamma^i_{jk}dx^jdx^k \psi_i-2dx^id\psi_i,
\end{equation}
where $\psi_1,\psi_2,\psi_3$ are an additional coordinates
(\cite{dryuma:paterson&walker})-(\cite{dryuma3:dryuma}).

    The Ricci tensor of the metric (\ref{dryuma:eq17}) is vanished
at the conditions (\ref{dryuma:eq2}-\ref{dryuma:eq7}), $\lambda=0$, and the Riemann
tensor $^{6}R_{i j k l}=0$ in this case also is vanished.

    So after the Rieman extension of the flat 3-dimensional diagonal metric $\lambda=0$ with help of
    the coefficients (\ref{dryuma:eq16}) we get a six-dimensional
   flat $^{6}R_{ijkl}=0$ space having the signature $[+++---]$.

 The Ricci tensor of the metric (\ref{dryuma:eq1}) at the
 conditions (\ref{dryuma:eq2}-\ref{dryuma:eq7})
has the components
\[
R_{11}=2\lambda A(x,y,z)^2, \quad R_{22}=2\lambda B(x,y,z)^2,\quad
R_{33}=2\lambda C(x,y,z)^2,\]\[\quad R_{12}=0,\quad R_{13}=0,\quad
R_{23}=0
\]
and corresponding space
   in the case $\lambda=\pm1$ is symmetric space with conditions
   on the curvature tensor
\[
^{3}R_{i j kl;m}=0.
\]

    It is remarkable fact that after the Riemann extension of the space of constant curvature
    one get symmetric space.

\begin{theor}
    The Riemann extension of 3-dimensional  space of constant curvature is a six-dimensional
     symmetric space
\[
^{6}R_{i j k l;m}=0.
\]
\end{theor}

 The proof of this statement can be checked by direct calculations.

\begin{rem}
  In the  case
\[
\Gamma^1_{23}\neq0,\quad\Gamma^3_{12}\neq0,\quad\Gamma^2_{13}\neq0
\]
it is possible to  get a six-dimensional curved metric
\begin{equation}\label{dryuma:eq18}
^{6}ds^2=-2\Gamma^i_{jk}dx^jdx^k \psi_i-2dx^id\psi_i,
\end{equation}
with conditions
\[
\Gamma^2_{13}=\frac {C(x,y,z)A(x,y,z)}{B(x,y,z)},\quad
\Gamma^3_{12}=\frac {A(x,y,z)B(x,y,z)}{C(x,y,z)},\quad
\Gamma^1_{23}=\frac {B(x,y,z)C(x,y,z)}{A(x,y,z)}
\]
for the spaces with $\lambda=1$
 and
\[
\Gamma^2_{13}=\frac {iC(x,y,z)A(x,y,z)}{B(x,y,z)},\quad
\Gamma^3_{12}=\frac {iA(x,y,z)B(x,y,z)}{C(x,y,z)},\quad
\Gamma^1_{23}=\frac {iB(x,y,z)C(x,y,z)}{A(x,y,z)}
\]
for the spaces with $\lambda=-1$.

  For such form of connection coefficients we get the generalization of
  the Darboux system of equations

\[
{\frac {\partial ^{2}}{\partial y\partial z}}A(x,y,z)={\frac
{\left ({ \frac {\partial }{\partial z}}B(x,y,z)\right ){\frac
{\partial }{
\partial y}}A(x,y,z)}{B(x,y,z)}}+{\frac {\left ({\frac {\partial }{
\partial y}}C(x,y,z)\right ){\frac {\partial }{\partial z}}A(x,y,z)}{C
(x,y,z)}}+\]\[+2\frac{\partial\left(B(x,y,z)C(x,y,z)\right)}{\partial
x} ,
\]
\begin{equation} \label{dryuma:eq19}
{\frac {\partial ^{2}}{\partial x\partial z}}B(x,y,z)={\frac
{\left ({ \frac {\partial }{\partial z}}A(x,y,z)\right ){\frac
{\partial }{
\partial x}}B(x,y,z)}{A(x,y,z)}}+{\frac {\left ({\frac {\partial }{
\partial x}}C(x,y,z)\right ){\frac {\partial }{\partial z}}B(x,y,z)}{C
(x,y,z)}}+\]\[+2\frac{\partial\left(A(x,y,z)C(x,y,z)\right)}{\partial
y} ,
\end{equation}
\[
{\frac {\partial ^{2}}{\partial x\partial y}}C(x,y,z)={\frac
{\left ({ \frac {\partial }{\partial y}}A(x,y,z)\right ){\frac
{\partial }{
\partial x}}C(x,y,z)}{A(x,y,z)}}+{\frac {\left ({\frac {\partial }{
\partial x}}B(x,y,z)\right ){\frac {\partial }{\partial y}}C(x,y,z)}{B
(x,y,z)}}+\]\[+2\frac{\partial\left(A(x,y,z)B(x,y,z)\right)}{\partial
z},
\]
which can be useful in applications.

    For example the simplest solution of this system is
\[
A(x,y,z)=\exp(z/4-x),\quad B(x,y,z)=\exp(y-z/4),\quad
C(x,y,z)=\exp(z/4)
\]
which corresponds the metrics (\ref{dryuma:eq18}) of a
six-dimensional manifold with a corresponding conditions on the
curvature tensor.

\end{rem}

\section{The normal Riemann spaces}

      The notion of normal Riemann space was introduced first by
      Eisenhart. Their properties have been studied in (\cite{dryuma:sob}).

\begin{dfn} An n-dimensional Riemannian space $N(x^k)$ with the metrics
\[
ds^2=g_{i j}dx^i dx^j
\]
is normal if a following conditions on their main curvatures $K_i$
hold
\[
\frac{\partial K_l}{\partial x^l}=3\lambda_l+3\mu_l K_l, \quad
i=l,
\]
and
\[
\frac{\partial K_i}{\partial
 x^l}=\lambda_l+\mu_l K_i,\quad i \neq l,
\]
\[
\frac{\partial \ln(g_{ij})}{\partial
x^l}=\frac{2}{K_l-K_i}\frac{\partial K_i}{\partial x^l} \quad i
\neq l,
\]
where $ \lambda$ and $\mu$ are the functions of coordinates $x^k$.
\end{dfn}

\begin{rem}

  The values $K_i$ are the roots of algebraic equation
\[
|b_{ij}-\kappa g_{ij}|=0
\]
relatively of a second range tensor $b_{ij}$ defined in a given
Riemann space.
\end{rem}

   In 3-dimensional case principal curvatures $K_i$ satisfy the system of equations
\[
(K_2-K_3)\frac{\partial K_1}{\partial x}+3(K_3-K_1)\frac{\partial
K_2}{\partial x}+3(K_1-K_2)\frac{\partial K_3}{\partial x}=0,
\]
\begin{equation} \label{dryuma:eq20}
3(K_2-K_3)\frac{\partial K_1}{\partial y}+(K_3-K_1)\frac{\partial
K_2}{\partial y}+3(K_1-K_2)\frac{\partial K_3}{\partial y}=0,
\end{equation}
\[
3(K_2-K_3)\frac{\partial K_1}{\partial z}+3(K_3-K_1)\frac{\partial
K_2}{\partial z}+3(K_1-K_2)\frac{\partial K_3}{\partial z}=0.
\]

      The relations between the systems of Darboux equations
      (\ref{dryuma:eq2})-(\ref{dryuma:eq4}) and the system
      (\ref{dryuma:eq20}) were established first in
      (\cite{dryuma:lame}).
\begin{theor}

The system of equations (\ref{dryuma:eq20}) for principal curvatures of
three-dimensional normal spaces transforms into the system of equations
\[
{\frac {\partial ^{2}}{\partial x\partial y}}{\it
K_1}(x,y,z)+{\frac {\left ({\frac {\partial }{\partial
y}}A(x,y,z)\right ){\frac {
\partial }{\partial x}}{\it K_1}(x,y,z)}{A(x,y,z)}}+\]\[+\left ({\frac {{
\frac {\partial }{\partial x}}A(x,y,z)}{A(x,y,z)}}+{\frac {{\frac
{
\partial }{\partial x}}B(x,y,z)}{B(x,y,z)}}-{\frac {{\frac {\partial ^
{2}}{\partial x\partial y}}A(x,y,z)}{{\frac {\partial }{\partial
y}}A( x,y,z)}}\right ){\frac {\partial }{\partial y}}{\it
K_1}(x,y,z)=0,
\]
\begin{equation} \label{dryuma:eq21}
{\frac {\partial ^{2}}{\partial x\partial z}}{\it
K_1}(x,y,z)+{\frac {\left ({\frac {\partial }{\partial
z}}A(x,y,z)\right ){\frac {
\partial }{\partial x}}{\it K_1}(x,y,z)}{A(x,y,z)}}+\]\[+\left ({\frac {{
\frac {\partial }{\partial x}}C(x,y,z)}{C(x,y,z)}}+{\frac {{\frac
{
\partial }{\partial x}}A(x,y,z)}{A(x,y,z)}}-{\frac {{\frac {\partial ^
{2}}{\partial x\partial z}}A(x,y,z)}{{\frac {\partial }{\partial
z}}A( x,y,z)}}\right ){\frac {\partial }{\partial z}}{\it
K_1}(x,y,z)=0,
 \end{equation}
\[
{\frac {\partial ^{2}}{\partial y\partial z}}{\it
K_1}(x,y,z)+\left ( {\frac {{\frac {\partial }{\partial
z}}B(x,y,z)}{B(x,y,z)}}+{\frac {{ \frac {\partial }{\partial
z}}A(x,y,z)}{A(x,y,z)}}-{\frac {{\frac {
\partial ^{2}}{\partial y\partial z}}A(x,y,z)}{{\frac {\partial }{
\partial y}}A(x,y,z)}}\right ){\frac {\partial }{\partial y}}{\it K_1
}(x,y,z)+\]\[+\left (-{\frac {\left ({\frac {\partial }{\partial
y}}A(x,y,z) \right ){\frac {\partial }{\partial z}}B(x,y,z)}{\left
({\frac {
\partial }{\partial z}}A(x,y,z)\right )B(x,y,z)}}+{\frac {{\frac {
\partial }{\partial y}}A(x,y,z)}{A(x,y,z)}}\right ){\frac {\partial }{
\partial z}}{\it K_1}(x,y,z)=0,
\]
if the following relations are hold
\[
{\frac {\partial }{\partial y}}{\it K_1}(x,y,z)-{\frac {\left ({
\frac {\partial }{\partial y}}A(x,y,z)\right )\left ({\it
K_2}(x,y,z)-{ \it K_1}(x,y,z)\right )}{A(x,y,z)}}=0,
\]
\[
{\frac {\partial }{\partial z}}{\it K_1}(x,y,z)+{\frac {\left ({
\frac {\partial }{\partial z}}A(x,y,z)\right )\left (-{\it
K_3}(x,y,z)+ {\it K_1}(x,y,z)\right )}{A(x,y,z)}}=0,
\]
\[
{\frac {\partial }{\partial x}}{\it K_2}(x,y,z)+{\frac {\left ({
\frac {\partial }{\partial x}}B(x,y,z)\right )\left ({\it
K_2}(x,y,z)-{\it K_1}(x,y,z)\right )}{B(x,y,z)}}=0,
\]
\[
{\frac {\partial }{\partial z}}{\it K_2}(x,y,z)+{\frac {\left ({
\frac {\partial }{\partial z}}B(x,y,z)\right )\left ({\it
K_2}(x,y,z)-{ \it K_3}(x,y,z)\right )}{B(x,y,z)}}=0,
\]
\[
{\frac {\partial }{\partial x}}{\it K_3}(x,y,z)-{\frac {\left ({
\frac {\partial }{\partial x}}C(x,y,z)\right )\left (-{\it
K_3}(x,y,z)+ {\it K_1}(x,y,z)\right )}{C(x,y,z)}}=0
\]
\[
{\frac {\partial }{\partial y}}{\it K_3}(x,y,z)-{\frac {\left ({
\frac {\partial }{\partial y}}C(x,y,z)\right )\left ({\it
K_2}(x,y,z)-{ \it K_3}(x,y,z)\right )}{C(x,y,z)}}=0.
\]
\end{theor}

\begin{rem}

   Similar equations and for the components $K_2$ and $K_3$
can be written.

   As example for the value $K_2(x,y,z)$ one get
\[
{\frac {\partial ^{2}}{\partial x\partial y}}{\it
K_2}(x,y,z)+{\frac {\left ({\frac {\partial }{\partial
x}}B(x,y,z)\right ){\frac {
\partial }{\partial y}}{\it K_2}(x,y,z)}{B(x,y,z)}}+\]\[+\left ({\frac {{
\frac {\partial }{\partial y}}B(x,y,z)}{B(x,y,z)}}+{\frac {{\frac
{
\partial }{\partial y}}A(x,y,z)}{A(x,y,z)}}-{\frac {{\frac {\partial ^
{2}}{\partial x\partial y}}B(x,y,z)}{{\frac {\partial }{\partial
y}}B( x,y,z)}}\right ){\frac {\partial }{\partial x}}{\it
K_2}(x,y,z)=0,
\]
\begin{equation} \label{dryuma:eq22}
{\frac {\partial ^{2}}{\partial x\partial z}}{\it
K_2}(x,y,z)+\left ({ \frac {{\frac {\partial }{\partial
z}}B(x,y,z)}{B(x,y,z)}}-{\frac { \left ({\frac {\partial
}{\partial x}}C(x,y,z)\right ){\frac {
\partial }{\partial z}}B(x,y,z)}{C(x,y,z){\frac {\partial }{\partial x
}}B(x,y,z)}}\right ){\frac {\partial }{\partial x}}{\it
K_2}(x,y,z)+ \]\[+\left ({\frac {{\frac {\partial }{\partial
x}}B(x,y,z)}{B(x,y,z)}}-{ \frac {{\frac {\partial ^{2}}{\partial
x\partial z}}B(x,y,z)}{{\frac {
\partial }{\partial z}}B(x,y,z)}}+{\frac {{\frac {\partial }{\partial
x}}C(x,y,z)}{C(x,y,z)}}\right ){\frac {\partial }{\partial z}}{\it
K_2} (x,y,z) =0,
\end{equation}
\[
{\frac {\partial ^{2}}{\partial y\partial z}}{\it
K_1}(x,y,z)+{\frac {\left ({\frac {\partial }{\partial
z}}B(x,y,z)\right ){\frac {
\partial }{\partial z}}{\it K_2}(x,y,z)}{B(x,y,z)}}+\]\[+\left ({\frac {{
\frac {\partial }{\partial y}}C(x,y,z)}{C(x,y,z)}}+{\frac {{\frac
{
\partial }{\partial y}}B(x,y,z)}{B(x,y,z)}}-{\frac {{\frac {\partial ^
{2}}{\partial y\partial z}}B(x,y,z)}{{\frac {\partial }{\partial
z}}B( x,y,z)}}\right ){\frac {\partial }{\partial z}}{\it
K_2}(x,y,z)=0.
\]
\end{rem}

    It is important to note that the equations (\ref{dryuma:eq21}) and (\ref{dryuma:eq22}) are present
    the linear systems of equations  with conditions of compatibility in form of the Darboux system
     (\ref{dryuma:eq2})-(\ref{dryuma:eq4}).

     This property will be used for the studying  of solutions of the complete Lame system of equations
     (\ref{dryuma:eq2})-(\ref{dryuma:eq7}).

   Now we present some solutions of the system
   (\ref{dryuma:eq20}).

     Remark that from the system (\ref{dryuma:eq20}) follows the
     relations
\[
K_i(x,y,z)(\phi_j-\phi_l)+K_j(x,y,z)(\phi_l-\phi_i)+K_l(x,y,z)(\phi_i-\phi_j)=0,\quad
i\neq j\neq l
\]
where $\phi_i=\phi_i(x^i)$ are arbitrary functions depending from the variable
$x^i$.

    This restrictive clause together with equations  (\ref{dryuma:eq20}) lead to the system of
    equations on the functions $K_i(x,y,z)$.

    As example for the value $K_1(x,y,z)$ we get
\[
{\frac {\partial ^{2}}{\partial x\partial y}}{\it
K_1}(x,y,z)-1/2\,{ \frac {-\left ({\frac {\partial }{\partial
x}}{\it K_1}(x,y,z)\right ){ \frac {d}{dy}}b(y)+3\,\left ({\frac
{d}{dx}}a(x)\right ){\frac {
\partial }{\partial y}}{\it K_1}(x,y,z)}{a(x)-b(y)}}=0,
\]
\[
{\frac {\partial ^{2}}{\partial x\partial z}}{\it
K_1}(x,y,z)-1/2\,{ \frac {3\,\left ({\frac {d}{dx}}a(x)\right
){\frac {\partial }{
\partial z}}{\it K_1}(x,y,z)-\left ({\frac {d}{dz}}c(z)\right ){\frac {
\partial }{\partial x}}{\it K_1}(x,y,z)}{a(x)-c(z)}}=0,
\]
\[
{\frac {\partial ^{2}}{\partial y\partial z}}{\it
K_1}(x,y,z)-1/2\,{ \frac {-\left ({\frac {\partial }{\partial
y}}{\it K_1}(x,y,z)\right ){ \frac {d}{dz}}c(z)+\left ({\frac
{\partial }{\partial z}}{\it K_1}(x,y, z)\right ){\frac
{d}{dy}}b(y)}{b(y)-c(z)}}=0
\]
where $a(x),\quad b(y),\quad c(z)$ are arbitrary functions.

  With the help of equations (\ref{dryuma:eq21}) it is possible to show that the components of metrics
   of the Darboux  space with a given conditions on principal curvatures are
\[
A(x,y,z)={\frac {U(x)}{\sqrt {a(x)-b(y)}\sqrt {a(x)-c(z)}}},
\]
\[
B(x,y,z)={\frac {V(y)}{\left (a(x)-b(y)\right )^{5/2}\sqrt
{b(y)-c(z)} }},
\]
\[
C(x,y,z)={\frac {W(z)}{\left (a(x)-c(z)\right )^{5/2}\sqrt
{b(y)-c(z)} }}.
\]


\begin{thebibliography}
\footnotesize
 \bibitem{dryuma:darboux} V.S.Dryuma, Projective properties of a family operators, IX- Geometric Conference, Kishinev, 20-22 Septembre,1988,
  Theses of communications, pp.104-105 (in Russian).

\bibitem{dryuma:darboux1} V.S.Dryuma,
 Three dimensional exactly integrable system of nonlinear equations and
 its applications,
 {\it Matematicheskie issledovaniya}, Kishinev, Stiintsa, 1992, v.124, pp.56-68 (in Russian).

\bibitem{dryuma:darboux2} V.S.Dryuma,
 Geometrical properties of the multidimensional nonlinear differential equations,{\it Teoreticheskaya i Matematicheskaya
  Fizika}, v.99, No.2, 1993, pp.241-249.

 \bibitem{dryuma:zakh} Zakharov V., Description of the n-orthogonal curvalinear coordinate systems
 and Hamiltonian integrable systems of hydrodynamic type: 1. Integration of the Lame Equations,
  {\it Duke Math. Journal}, 1998, v.94 No.1,p.103--139.

\bibitem{dryuma:wolf} Wolf T.,  About vacuum solutions of Einstein's field equation with flat three-dimensional
hypersufaces,{\it  Journal of Math. Phys., v.27(9)}, 1986,
2354-2359.

 \bibitem{dryuma:jac} Jackiw R., A Pure Cotton Kink in a Funny Place,
 {\it ArXiv: math-ph/0403044}, v.2, 21 July 2004.

 \bibitem{dryuma:paterson&walker} Paterson E.M. and Walker A.G.,Riemann extensions,
 {\it Quart.J.Math.Oxford},1952, V.3,19--28.

 \bibitem{dryuma1:dryuma} Dryuma V.,The Riemann Extensions in theory of differential equations and their
  applications, {\it Matematicheskaya fizika, analiz, geometriya}, 2003,v.10, No.3,1--19.

 \bibitem{dryuma2:dryuma} Dryuma V., The Riemann and Einstein-Weyl geometries in the theory of ODE's
 , their applications and all that,
 {New Trends in Integrability and Partial Solvability ,115-156, (eds.A.B.Shabat et al.) 2004, Kluwer Academic Publishers},
{ ArXiv: nlin: SI/0303023, 11 March, 2003, 1--37}.

 \bibitem{dryuma3:dryuma} Dryuma V.,Applications of Riemannian and Einstein-Weyl Geometry in
 the theory of second order ordinary differential equations,{\it Theoretical and
 Mathematical Physics}, 2001, V.128, N 1, 845--855.

\bibitem{dryuma:sob} Sobchuk V.,  About one classe of normal Riemann spaces
,{\it Seminar on vector and tensor analyse}, v.15, 1968, 65-76 (in Russian).

\bibitem{dryuma:lame} Dryuma V., On the law of transformation of affine connection and its integration. Part 1.
 Generalization of the Lame equations, {\it
 Buletinul Academiei de Stiinte a Republicii Moldova, matematika}, 1998, V.1(26), 55--68.
 \end{thebibliography}
 \end{document}